\documentclass{agtart_a} 
\pdfoutput=1

\usepackage{graphicx}
\usepackage{pinlabel}

\title{From continua to $\mathbb{R}$--trees}

\author{Panos Papasoglu}
\givenname{Panos}
\surname{Papasoglu}
\address{Mathematics Department\\
University of Athens\\\newline
Athens 157 84\\Greece}
\email{panos@math.uoa.gr}
\urladdr{}

\author{Eric Swenson}
\givenname{Eric L}
\surname{Swenson}
\address{Mathematics Department\\
Brigham Young University\\\newline
Provo UT 84602\\USA}
\email{eric@math.byu.edu} 
\urladdr{}

\volumenumber{6}
\issuenumber{}
\publicationyear{2006}
\papernumber{62}
\startpage{1759}
\endpage{1784}

\doi{}
\MR{}
\Zbl{}


\received{31 May 2006}
\revised{25 August 2006}
\accepted{26 August 2006}
\published{1 November 2006}
\publishedonline{1 November 2006}
\proposed{}
\seconded{}
\corresponding{}
\editor{Liz}
\version{}

\arxivreference{} 

\AtBeginDocument{\let\bar\wbar\let\tilde\wtilde\let\hat\what}
\AtBeginDocument{\newcommand{\interior}{^{ \kern-5pt ^\circ}}}
\let\Bbb\mathbb
\makeop{fix}
\def\hatphi{\mskip2mu\hat{\mskip-2mu\vphantom{t}\smash{\phi}\mskip-.5mu}\mskip.5mu}

\makeatletter
\def\cnewtheorem#1[#2]#3{\newtheorem{#1}{#3}
\expandafter\let\csname c@#1\endcsname\c@Thm}
\makeatother

\newtheorem{Thm}{Theorem}
\makeautorefname{Thm}{Theorem}
\cnewtheorem{Cor}[Thm]{Corollary}
\makeautorefname{Cor}{Corollary}
\cnewtheorem{Lem}[Thm]{Lemma}
\makeautorefname{Lem}{Lemma}
\cnewtheorem{Prop}[Thm]{Proposition}

\theoremstyle{definition}
\newtheorem{Def}{Definition}

\theoremstyle{remark}

\cnewtheorem{Ex}[Thm]{Example}
\newtheorem{No}{Notation}

\newcommand {\bd}{\partial}

\newcommand {\iy}{\infty}
\newcommand {\N}{{\mathbb N}}

\newcommand {\diam}{\text{diam}}

\newcommand {\cB}{ {\mathcal  B}}

\newcommand {\cC}{{\mathcal  C}}

\newcommand {\cP}{{\mathcal  P}}
\newcommand {\cR}{{\mathcal  R}}

\begin{document}

\begin{asciiabstract}
We show how to associate an R-tree to the set of cut
points of a continuum. If X is a continuum without cut points we
show how to associate an R-tree to the set of cut
pairs of X.
\end{asciiabstract}

\begin{htmlabstract}
We show how to associate an <b>R</b>&ndash;tree to the set of cut
points of a continuum. If X is a continuum without cut points we
show how to associate an <b>R</b>&ndash;tree to the set of cut
pairs of X.
\end{htmlabstract}

\begin{abstract}
We show how to associate an $\mathbb{R} $--tree to the set of cut
points of a continuum. If $X$ is a continuum without cut points we
show how to associate an $\mathbb{R} $--tree to the set of cut
pairs of $X$.

\end{abstract}

\maketitle

\section {Introduction}

The study of the structure of cut points of continua has a long
history. Whyburn \cite{WHY} in 1928 showed that the set of cut
points of a Peano continuum has the structure of a ``dendrite''.
This ``dendritic'' decomposition of continua has been extended and
used to prove several results in continua theory.

We recall here that a continuum is a compact, connected metric
space and a Peano continuum is a locally connected continuum. If
$X$ is a continuum we say that a point $c$ is a cut point of $X$
if $X-\{c \}$ is not connected.

Continua theory became relevant for group theory after the
introduction of hyperbolic groups by Gromov \cite{GRO}. The
Cayley graph of a hyperbolic group $G$ can be ``compactified'' and
if $G$ is one-ended its Gromov boundary $\partial G$ is a
continuum. Moreover the group $G$ acts on $\partial G$ as a
convergence group. It turns out that algebraic properties of $G$
are reflected in topological properties of $\partial G$. A
fundamental contribution to the understanding of the relationship
between $\partial G$ and algebraic properties of $G$ was made by
Bowditch. In \cite {BOW1} Bowditch shows how to pass from the
action of a hyperbolic group $G$ on its boundary $\partial G$ to
an action on an $\mathbb{R}$--tree. The construction of the tree
(under the hypothesis of the $G$--action) from the continuum is
similar to the dendritic decomposition of Whyburn. The difficulty
here comes from the fact that the continuum is not assumed to be
locally connected.

The second author in \cite{SWE} (see also \cite{SWE3}) explained
how to associate to any continuum a ``regular big tree'' $T$, and
conjectured that $T$ is in fact an $\mathbb{R}$--tree. It is this
conjecture that we prove in the first part of this paper.

Let $X$ be a continuum without cut points. If $a,b\in X$ we say
that $a,b$ is a {\em cut pair\/} if $X-\{a,b \}$ is not connected. In
the second part of this paper we show how to associate an
$\mathbb{R}$--tree to the set of cut pairs of $X$ (compare \cite{BOW6}). We call this
tree a JSJ-tree motivated by the fact that if $G$ is a one-ended
hyperbolic group then the tree associated to $\partial G$ by this
construction is the tree of the JSJ-decomposition of $G$ (in this
case one obtains in fact a simplicial tree). Continua appear in
group theory also as boundaries of ${\rm CAT}(0)$ groups. In
\cite{P-S} we use the construction of $\mathbb{R}$--trees from cut
pairs presented here to extend Bowditch's results on splittings
\cite{BOW1} to ${\rm CAT}(0)$ groups. We show in particular that if
$G$ is a one-ended ${\rm CAT}(0)$ group such that $\partial G$ has a cut
pair then $G$
splits over a virtually cyclic group.

\subsubsection*{Acknowledgements} We would like to thank the referee for many suggestions that
improved the exposition and for correcting a mistake in the proof
of \fullref{T:countable}.

 This work is co-funded by European Social Fund and National
Resources (EPEAEK II) PYTHAGORAS.

\section {Preliminaries}
\begin{Def}[Pretrees]
Let $\cP$ be a set with a betweenness relation. If $y$ is between
$x$ and $z$ we write $xyz$. $\cP$ is called a pretree if the following
hold:

1. There is no $y$ such that $xyx$ for any $x\in \cP$.

2. $xzy \Leftrightarrow yzx$.

3. For all $x,y,z$, if $y$ is between $x$ and $z$ then $z$ is not between
$x$ and $y$.

4. If $xzy$ and $z\ne w$ then either $xzw$ or $yzw$.
\end{Def}

\begin{Def} We say that a pretree $\cP$ is discrete if for any
$x,y\in \cP$ there are finitely many $z\in \cP$ such that $xzy$.
\end{Def}

\begin{Def}A compact connected metric space is called a continuum.
\end{Def}
\begin{Def} Let $X$ be a topological space. We say
that a set  $C$ \emph{separates\/} the nonempty sets  $A,B \subset X $
if there are disjoint open sets $U,V$ of $X-C$, such that
$A\subset U$, $B \subset V$ and $U \cup V =X- C$. We say $C$
separates the points $a,b \in X$ if $C$ separates $\{a\}$ and
$\{b\}$. We say that $C$ separates $X$ if $C$ separates two points
of $X$. If $C=\{c\}$ then we call $c$ a cut point.  If $C=
\{c,d\}$ where $c\neq d$ and \emph{neither\/} $c$ nor $d$ is cut
point, then we call $\{c,d\}$ a (unordered) cut pair.
\end{Def}
The proof of the following Lemma is an  elementary exercise in
topology and will be left to the reader.

\begin{Lem} \label{L:etop}
Let $A$ be a connected subset of the space $X$ and $B$ closed in $X$.  If $A
\cap \mathrm{Int} B \neq \emptyset$,
then either $A \subset B$ or $A \cap \bd B$ separates the subspace $A$.
\end{Lem}

\begin{Lem} \label{L:eas} Let $X$ be a continuum and $C \subset X$ be
minimal with the property that $X-C$ is not connected.
The set   $C$ separates $A \subset X-C$ from $B \subset  X-C$ if and only if
there exist continua
$Y,Z$ such that
$A \subset Y$, $B \subset Z$,  $Y \cup Z = X$ and $Y \cap Z = C$.
\end{Lem}
\begin{proof}
We first show that $C$ is closed in $X$.   There are disjoint nonempty
subsets $D$ and $E$ open in $X-C$ with $D \cup E = X-C$.  By symmetry, it
suffices to show that $D$ is open in $X$. Suppose that $d \in D \cap \bd  D$.
There is a neighborhood $G$ of $d$ in $X$ such that $\bar G \cap \bar E =
\emptyset$.
Since $G \not \subset D$,  there is $c \in C \cap G$.  Notice that $D \cup
\{c\}$ and $E$ are disjoint open subsets of $X -(C-\{c\})$ with $(D\cup
\{c\}) \cup E = X-(C-\{c\})$, and $C$ is not minimal.
Therefore $C$ must be closed.

Suppose now that $C$ separates $A$ from $B$.  Thus there exist disjoint
nonempty  $U$ and $V$ open in $X-C$ (this implies open in $X$) such that $A
\subset U$, $B \subset V$, $U \cup V = X-C$.  Since $\bd U$ separates $X$ ,
by the minimality of $C$, $\bd U = C = \bd V$.
 Suppose the closure $\bar U$ is not connected.  Then $\bar U = P \cup Q$
where $P$ and $Q$ are disjoint nonempty clopen (closed and open) subsets of
$\bar U$.  Since $\bar U$ is closed, this implies that $P$ and $Q$ are
closed subsets of $X$.  Since $\bar U \not \subset C$, we may assume that
$P\not \subset C$.  The boundary of $P$ in $\bar U$ is empty, so $\bd P
\subset \bd \bar U = \bd U =C$.  Again by minimality $\bd P = C$.  Since $P$
is closed in $X$,  $C \subset P$.  Thus $Q \subset U$, and $Q$ is open in
$U$ since it is open $\bar U$.  Thus $Q$ is clopen in $X$ which contradicts
$X$ being connected.

 The implication in the other direction is trivial.
\end{proof}

This next result is just an application of the previous result.
\begin{Lem} \label{L:top}
Let $X$ be a continuum  and $A,B \subset X$.
\begin{itemize} \item The point $c \in X$ is
a cut point of $X$ which separates $A$ from $B $ if and only if
there exist continua $Y,Z \subset X$ such that $A \subset
Y-\{c\}$, $B \subset Z-\{c\}$, $Y\cup Z= X$ and $Y \cap Z =
\{c\}$. \item The pair of non-cut points $\{c,d\}$ is a cut pair
separating $A$ from $B $ if and only if there exist continua $Y,Z
\subset X$ such that $A \subset Y-\{c,d\}$, $B\subset Z-\{c,d\}$,
$Y\cup Z = X$, $Y \cap Z = \{c,d\}$.
\end{itemize}
\end{Lem}
\section{Cutpoint trees} Let $X$ be a metric continuum.
In this section we show that the big tree constructed in \cite{SWE3}
is always a real tree.
For the reader's convenience we recall briefly the construction
here.

For the remainder of this section,  $X$ will be a continuum.
\begin{Def}
If $a, b \in X$ we say that $c\in (a,b)$ if the cut point $c$ separates $a$ 
from $b$.

We call $(a,b)$ an interval and this relation an interval
relation. We define closed and half open intervals in the obvious
way, ie $[a,b)=\{a \}\cup (a,b)$, $[a,b]=\{a,b \}\cup (a,b)$ for
$a\ne b$ and $[a,a)=\emptyset, \, [a,a]=\{a \} $.
\end{Def}
\begin{Def} We define an equivalence relation on $X$. Each cut point
is equivalent only to itself and if
$a, b\in X$ are not cut points we say that $a$ is equivalent to $b$
 if
$(a,b)=\emptyset $.
\end{Def}
Let's denote by $\cP$ the set of equivalence classes for this
relation. We can define an interval relation on $\cP$ as follows:
\begin{Def} If $x, y \in \cP$ and $c$ is a cut point (so $c \in \cP$)
we say that $c \in (x,y)$ if for some (any) $a\in x$, $b \in y$,  we have
$c \in (a,b)$.

For $z\in \cP$, $z$ not a cut point we say that $z\in (x,y)$ if for some 
(any) $a\in
x,b\in y, c\in z$ we have that
$$[a,c)\cap (c,b]=\emptyset. $$
\end{Def}

If $x, y,z\in \cP$ we say that $z$ is between $x,y$ if $z\in
(x,y)$. We will show that $\cP$ with this betweenness relation is a
pretree. The first two axioms of the definition of pretree are
satisfied by definition.

For the remaining two axioms we recall the following lemmas (for a
proof see Bow\-ditch \cite{BOW5} or Swenson \cite{SWE}).

\begin{Lem} \label{L:ax3}
For any $x,y \in \cP$, if $z\in (x,y)$ then $x\notin (y,z)$.
\end{Lem}
\begin{Lem} \label{L:ax4}
For any $x,y,z\in P$, $(x,z)\subset (x,y]\cup [y,z)$.
\end{Lem}
Axiom 3 follows from \fullref{L:ax3} and Axiom 4 from \fullref{L:ax4}.

Now consider the following example where $X\subset \R^2$ is the
union  of a Topologist's sine curve, two arcs, five circles and
two disks:

\begin{center}
\labellist
\pinlabel $X$ at 256 43
\endlabellist
\includegraphics[scale=1.2]{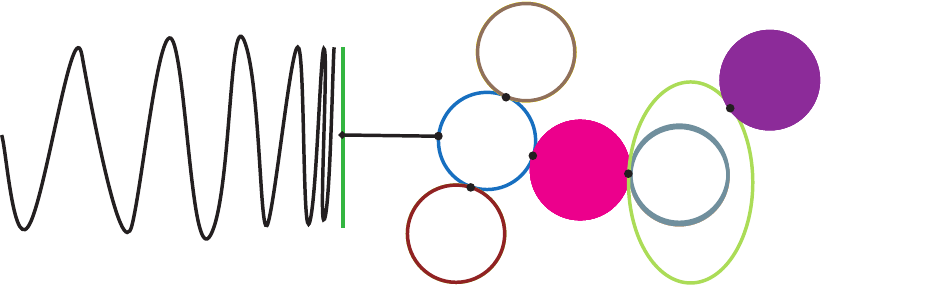}

\labellist 
\pinlabel $\cP$ at 250 32
\endlabellist
\includegraphics[scale=1.2]{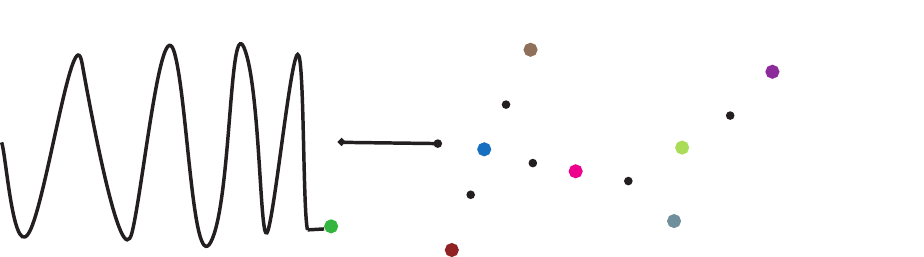}

\labellist 
\pinlabel $T$ at 249 32
\endlabellist
\includegraphics[scale=1.2]{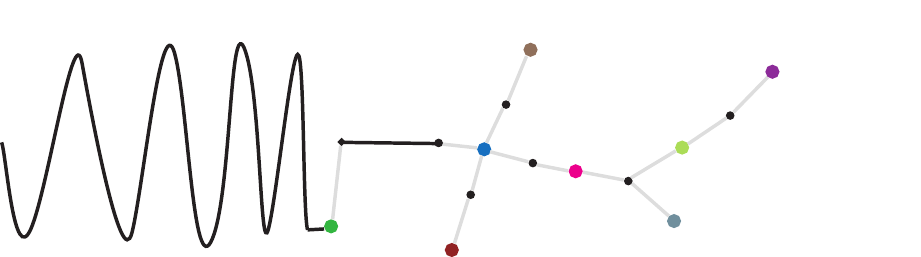}
\end{center}

The tree $T$ is obtained from $\cP$ by ``connecting the dots" according to
the pretree relation on $\cP$.
We will give the rigorous definition of $T$ later.

We have the following results about intervals in pretrees from
\cite{BOW5}:

\begin{Lem} \label{L:subset}
If $x,y,z \in \cP$, with $y\in [x,z]$ then $[x,y]\subset [x,z]$.
\end{Lem}

\begin{Lem} \label{L:order}
Let $[x,y]$ be an interval of $\cP$. The interval structure
induces two linear orderings on $[x,y]$, one being the opposite of
the other, with the property that if $<$ is one of the orderings,
then for any $z$ and $w\in [x,y]$ satisfying $z<w$, we have $(z,w)=\{ u\in [x,y]:z<u<w
\}$. In other words the interval structure defined by one of the
orderings is the same as our original interval structure.
\end{Lem}
\begin{Def} If $x$ and $y$ are distinct points of $\cP$ we say that $x$ and $y$
are \textit{adjacent\/} if $(x,y)=\emptyset $.
We say $x \in \cP$ is {\em terminal\/} if there is no pair $y,z \in \cP$ with
$x \in (y,z)$.
\end{Def}
We recall the following lemma from \cite{SWE}.
\begin{Lem} \label{L:adjacent}
If $x,y \in \cP$, are adjacent then exactly one of them is a cut
point and the other is a nonsingleton equivalence class whose
closure contains this cut point.
\end{Lem}
\begin{Cor}\label{L:singing} If $p \in \cP$ is a singleton equivalence class
and
$p$ is not a cut point, then $p$ is terminal in $\cP$.
\end{Cor}
\begin{proof}
Let $x \in X$ with $[x] =\{x\}$, and $x \in (a,b)$ for some $a,b \in \cP$.
Suppose that $x$ is not a cut point.  By \fullref{L:adjacent} there is no
point of $\cP$ adjacent to $[x]$.  Thus there are infinitely many cut points 
in $[a,x]$.
For each such cut point $c \in (a,x) $ choose a  continuum
$A_c \ni a$ with $\bd A_c = \{c\}$.
Considered the nested union $A = \bigcup A_c$.
We will show that $\bd A = \{x\}$.

  First consider $ y \in A$.  There exists  a cut point $c\in (a,x)$, $c 
\neq y$,  with
$y \in A_c$.  Thus $y \in \mathrm{Int}
A_c$ so
  $y \not \in \bd A$.

  Now consider $z \in X-A$ with $z \neq x$.   Since $z \not \in A$, by
definition $x \in ([z], [y])$ for any $y \in A$.  Since $x$ and $[z]$ are 
not adjacent there is a cut point $d \in ([z],x) $.
There exist continua $Z, B$ such that $Z \cup B = X$, $z \in Z$, $x \in
B$ and $Z\cap B =\{d\}$.  Since $x \in ([z], [y])$ for any $y \in A$, by
definition $A \subset B$, and $z \not \in \bd A$.

  The fact that $x\in \bd A$ follows  since $b \not \subset A$, so  $X \neq
A$, and so $\bd A \neq \emptyset$.
\end{proof}
We have the theorem \cite[Theorem 6]{SWE}:
\begin{Thm}\label{T:nested}
A nested union of intervals of $\cP$ is an interval of $\cP$.
\end{Thm}
\begin{Cor}\label{C:supremum} Any interval of $\cP$ has the supremum
property with respect to either of the linear orderings derived
from the interval structure.
\end{Cor}
\begin{proof}
Let $[x,y]$ be an interval of $\cP$ with the linear order $\leq $.
Let $A\subset [x,y]$. The set $\{[x,a]:a \in A\}$ is a set of nested
intervals so their union is an interval $[x,s]$ or $[x,s)$ and
$s=\sup\,A$.
\end{proof}
\begin{Def} A {\em big arc\/} is the homeomorphic image of a compact connected
nonsingleton linearly ordered topological space.  A separable big arc is
called an {\em arc\/}. A {\em big tree\/} is a  uniquely big-arcwise connected
topological space.  If all the big arcs of a big tree are arcs, then the big
tree is called a {\em real tree\/}.   A metrizable real tree is called an
$\R$--tree. An example of a real tree which is not an
$\mathbb{R}$--tree is the long line \cite[Section 2.5, page 56]{HOC-YOU}.
\end{Def}

\begin{Def} A pretree $\cR$ is {\em complete\/} if  every closed interval  is
complete as a
linearly ordered topological space (this is slightly weaker than the
definition given in \cite{BOW5}). Recall that a linearly ordered topological
space is complete if every bounded set has a supremum.

Let $\cR$ be a pretree.
An interval $I\subset \cR$ is called {\em preseparable\/} if
there is a countable set $Q \subset I$ such that for every nonsingleton
closed interval $J \subset I$, we have
$J \cap Q \neq \emptyset$.
A pretree is {\em preseparable\/} if every interval in it is preseparable.
\end{Def}
We now give a slight generalization of a construction in \cite{SWE,SWE3}.  Let $\cR$ be a complete pretree.   Set
$$T=\cR \cup \underset {{x,\,y\, \text{adjacent}}} { \bigsqcup } I_{x,y}$$
where $I_{x,y}$ is a copy of the real open interval $(0,1)$ glued
in between $x$ and $y$. We extend the interval relation of $\cR$
to $T$ in the obvious way (as in \cite{SWE,SWE3}), so that
in $T$, $(x,y) = I_{x,y}$.    It is clear that $T$ is a complete
pretree with no adjacent elements.   When $\cR=\cP$, we call the
$T$ so constructed the cut point tree of $X$.

\begin{Def} For $A$ finite subset of $T$ and $s\in T$ we define
$$U(s,A)=\{ t\in T:[s,t]\cap A=\emptyset \}.$$
\end{Def}

The following is what the proof of \cite[Theorem 7]{SWE}   proves in this
setting.
\begin{Thm} Let $\cR$ be a complete pretree.  The pretree $T$, defined
above, with the topology defined by the basis  $\{U(s,A)\}$ is a regular big
tree.  If in addition $\cR$ is preseparable, then
$T$ is a real tree.
\end{Thm}

We now prove the conjecture from \cite{SWE}.
\begin{Thm} \label{T:countable}
The pretree $\cP$ is preseparable, so the cut point tree $T$ of
$X$ is a real tree.
\end{Thm}
\begin{proof}  By the proof of \cite[Theorem 7]{SWE}, it suffices to show
that there are only countably many adjacent pairs in a closed
interval $[a,b]$ of $\cP$. By \fullref{L:adjacent}, adjacent
elements of $\cP$ are pairs $E,c$ where $E$ is a nonsingleton
equivalence class, $c$ is a cut point and $c \in \bar E-E$. Let's
assume that there are uncountably many such pairs in $[a,b]$.
By symmetry we may assume that $E \in (a,c)$ for uncountably many pairs 
$(E,c)$, and for
each such pair we pick an $e\in E$.

Since $c$ separates $e$ from $b$ choose continua $A$, $B$ such that
$X = A \cup B$, $\{c\} = A \cap B$,  $e \in A$ and $b \in B$.
Since $e \not \in B$ and $B$ is compact  $d(e,B) >0$. Let $\epsilon_e=d(e,B)$.

In this way to each pair $E,c$  we
associate $e\in E$  continua $A$, $B$ and $\epsilon_e >0$.
Since there are uncountably many $e$, for some $n \in N$
there are  uncountably many $e$ with $\epsilon_e>1/n$.
Let's denote by $S$ the set of all such $e$ with $\epsilon_e >1/n$. Consider 
a finite
covering of $X$ by open balls of radius $\unpfrac{1}{2n}$. Since $S$
is infinite there are distinct elements $e_1,e_2,e_3\in S$ lying
in the same ball. It follows that $d(e_i,e_j)<1/n$ for all $i,j$.

The points $e_1,e_2,e_3$ correspond to adjacent elements of $\cP$, say
$E_1,c_1$,
$E_2,c_2,E_3,c_3$.  Since all these lie in an interval of
$\cP$ they are linearly ordered and we may assume, without
loss of generality, that $E_1\in [a,c_2)$ and $E_3\in
(c_2,b]$.  Let $A_1$ and $B_1$ be the continua chosen for $E_1,c_1$  such 
that $A_1 \cap B_1 = \{c_1\}$, $A_1\cup B_1 = X$,  $e_1 \in A_1$, $b 
\in B_1$ and $d(e_1,B_1) =\epsilon_{e_1}>1/n$.
 It follows that $e_3 \in B_1$ and so $d(e_1,e_3)\ge d(e_1,B_1) > 1/n$, which is a 
contradiction.
\end{proof}

The real tree $T$ is not always metrizable.  Take for example $X$
to be the cone on a Cantor set $C$ (the so-called Cantor
fan). Then $X$ has only one cut point, the cone point $p$, and
$\cP $
has uncountable many other elements $q_c$, one for each point $c \in C$.
As a pretree,
$T$ consists of uncountable many arcs $\{ [p,q_c]:\, c \in C\}$ radiating
from a single
central point $p$.  However, in the topology defined from the basis
$\{U(s,A)\}$, every
open set containing $p$ contains the arc $[p,q_c]$ for all but finitely
many $c \in C$.
There can be no metric, $d$, giving this topology since $d(p,q_c)$ could
only be nonzero
for countably many $c \in C$.

It is possible however to equip $T$ with a metric that preserves
the pretree structure of $T$. This metric is ``canonical'' in the
sense that any homeomorphism of $X$ induces a homeomorphism of
$T$. The idea is to metrize $T$ in two steps.
In the first step
one metrizes the subtree obtained by the span of cut points of
$\cP$. This can be written as a countable union of intervals and
it is easy to equip with a metric.

$T$ is obtained from this tree
by gluing intervals to some points of $T$. In this step one might
glue uncountably many intervals but the situation is similar to
the Cantor fan described above. The new intervals are metrized in
the obvious way, eg one can give all of them length one.

\begin{Thm}\label{T:Rtree}
There is a path metric $d$ on $T$, which preserves the pretree
structure of $T$, such that $(T,d)$ is a metric $\mathbb{R}
$--tree. The topology so defined on $T$ is  canonical (and may be
different from the topology with basis $\{U(s,A)\}$). Any
homeomorphism $\phi $ of $X$ induces a homeomorphism $\hatphi$
of $T$ equipped with this metric.
\end{Thm}
\begin{proof}
Let $\cC$ be the set of cut points of $X$ and let $S$ be a
countable dense subset of $\cC$.
 Choose a base point $s\in S$. Denote by $T'$ the
union of all intervals $[s,s']$ of $T$ with $s' \in S$. Now we
remark that at most countably many cut points of $X$ are not
contained in $T'$. Indeed if $c\in \cC$ is a cut point not in $T'$
then $X-c=U\cup V$ where $U,V$ are disjoint open sets and one of
the two (say $U$) contains no cut points. Let $\epsilon >0$ be
such that a ball $B(c)$ in $X$ of radius $\epsilon $ is contained
in $U$. So we associate to each $c$ not in $T'$ a ball $B(c)$ and
we remark that to distinct $c$'s correspond disjoint balls.
Clearly there can be at most countably many such disjoint balls in
$X$. Thus by enlarging $S$ we may assume that $T'$ contains all
cut points, so $T'$ is the convex hull of $\cC$ in $T$, and so is
canonical.

Since $S$ is countable we can write $S=\{s_1,s_2, \dots \}$ and we
metrize $T'$ by an inductive procedure: we give $[s,s_1]$ length 1
(Choose $f\co [0,1]\to [s,s_1]$, a homeomorphism, and define
$d(f(a),f(b))=|a-b|$). $[s,s_2]$ intersects $[s,s_1]$ along a
closed interval $[s,a]$. If $[a,s_2]$ is non empty we give it
length $1/2$ and we obtain a finite tree. At the $n$--th step of the
procedure we add the interval $[s,s_{n+1}]$ to a finite tree
$T_n$. If $[s,s_{n+1}]\cap T_n=[a,s_{n+1}]$ a nondegenerate
interval we glue $[a,s_{n+1}]$ to $T_n$ and give it length
$1/2^n$. Note that if $a,b\in T'$ then $a\in [s,s_n], b\in
[s,s_k]$ for some $k,n\in \Bbb N$.  Without loss of generality,
$k\le n$ and so $a,b\in T_n$ and $d(a,b)$ is determined at some
finite stage of the above procedure.

We remark that each end of $ T'$ corresponds to an element of $\cP
$ and by adding these points to $ T'$ one obtains an $\Bbb
R$--tree that we still denote by $ T'$. Here by end of $T'$ we mean
an ascending union $\bigcup [s,s_i]\ (i\in \Bbb N,\,s_i\in S)$ which
is not contained in any interval of $ T'$. If $C_i$ is the closure
of the union of all components of $X-s_i$ not containing $s$
we have that $C_i\subset C_{i-1}$ for all $i$ and $\bigcap C_i$ is
an element of $\cP$.

If $x\in \cP$ does not lie in $T'$ then $x$ is adjacent to
some cut point $c\in T'$. For each such adjacent pair $(c,x)$, by
construction $(c,x)$ is a copy of the unit interval $(0,1)$ and this gives
the path metric on $[c,x]$. In this way we equip $T$ with a path metric
$d$.

Clearly a homeomorphism $\phi\co X\to X$   induces a pretree
isomorphism $\hatphi\co \cP \to \cP$. By extending $\hatphi$ to
the intervals corresponding to adjacent points of $\cP$ (via the
identity map on the unit interval, $(0,1) \to(0,1)$)  we get a
pretree isomorphism  function $\hatphi\co T \to T$ which restricts
to a pretree isomorphism $\hatphi\co T' \to T'$. For any (possibly
singleton) arc $\alpha \in T'$  let $\cB _{\alpha }$ be the set of
components of $T' -\alpha$. By the construction of
$d$, for any $\epsilon >0$, the set $\{B \in \cB _{\alpha }:\,
\diam(B) > \epsilon\}$ is finite. It follows that $\hatphi\co  T'
\to T'$ is continuous (using the metric $d$) and therefore a
homeomorphism. We extend $\hatphi$ to $T$ by defining it to  be
an isometry on the disjoint union of  intervals $T-T'$. Thus we
get $\hatphi\co T \to T$ a homeomorphism.
\end{proof}

%

\section{JSJ-trees}

\begin{Def} Let $X$ be a continuum without cut points.
A finite set $S\subset X$ with
$|S|>2$ is called {\em cyclic subset\/}
if  there is an ordering $S=\{x_1,\dots x_n\}$ and continua $M_1,\dots M_n$
with the following properties:
\begin{itemize}
\item  $ M_n \cap M_1 = \{x_1\} $,  and for $i>1$, $\{x_i\} = M_{i-1} \cap
M_{i}$
\item $M_i \cap M_j = \emptyset$ for $i-j \neq \pm 1 \mod n$
\item $\bigcup M_i = X$
\end{itemize}
The collection $M_1, \dots M_n$ is called the (a) {\em cyclic decomposition\/}
of $X$ by $\{x_1,\dots x_n\}$.  This decomposition is unique as we show in \fullref{L:cross} (for $n>3$).

We also define a cut pair to be {\em cyclic\/}.

Clearly  every nonempty nonsingleton subset of a cyclic set
is cyclic.

If $S$ is an infinite subset of $X$ and every finite subset $A \subset S$
with $|A|>1$ is cyclic,
then we say $S$ is {\em cyclic\/}.
\end{Def}

Clearly if $A$ is a subset of a cyclic set with $|A|>1$, then $A$ is cyclic.

\begin{Lem} \label{L:cross}
Let $X$ be a connected metric space without cut points.
  If the cut pair ${a,b}$ separates the
cut pair $c,d $ then $\{a,b,c,d\}$ is cyclic, so $\{c,d\}$
separates $a$ from $b$. Furthermore $X-\{c,d\}$ has exactly two
components and $X-\{a,b\}$ has exactly two components.
\end{Lem}
\begin{proof}
By \fullref{L:top} there exist continua $C,D$ with $C\cap D =
\{a,b\}$, $X= C \cup D$ , $c\in C$ and $d \in D$. Since $c,d$ is a
cut pair there exist continuum $A,B$ such that $A\cup B = X$ and
$A\cap B = \{c,d\}$.  We may assume that $a \in A$.  Since $B$ is
connected, with $c,d \in B$ and $a\not \in B$, then $b$ must be a
cut point of $B$ separating $c$ and $d$.  Similarly $a$ is a cut
point of $A$ separating $c$ and $d$.  Thus by \fullref{L:top}
there exist continua $M_{a,c}$, $M_{a,d}$ $M_{b,c}$, $M_{b,d}$
with $c \in M_{a,c}$, $c \in M_{b,c}$,  $d \in M_{a,d}$, and $d
\in M_{b,d}$, such that $M_{a,c}\cup M_{a,d} =A $, $M_{a,c}\cap
M_{a,d} =\{a\}$, $M_{b,c}\cup M_{b,d} =B $ and  $M_{b,c}\cap
M_{b,d} =\{b\} $.  It follows that $M_{a,c} \cap M_{b,c} =\{c\}$
and that $M_{a,d} \cap M_{b,d} =\{d\}$.  Thus $\{a,b,c,d\}$ is
cyclic.

Suppose that $\{c,d\}$ separated $A$, then there would be nonsingleton
continua $F$,$G$ with $A = F\cup G$ and $\{c,d\}= F\cap G$. We may assume
that $a \in F$. Since $a$ separates $c$ from $d$ in $A$, either $c \not \in
G$ or $d \not \in G$.  With no loss of generality $d \not \in G$.  Thus $F
\cup B$ and $G$ are continua with $X= (F\cup B) \cup G$ and $\{c\}= (F\cup
B) \cap G$, making $c$ a cut point of $X$.  This is a contradiction, so
$\{c,d\}$ doesn't separate $A$ and similarly $\{c,d\}$ doesn't separate $B$.
Thus $X-\{c,d\}$ has exactly two components.
\end{proof}
\begin{Def} Let $X$ be a metric space without cut points.  A
nondegenerate nonempty set $A \subset X$ is called inseparable if
no pair of points in $A$ can be separated by a cut pair.
\end{Def}
Every inseparable set is contained in a maximal inseparable set. A
maximal inseparable subset is closed (its complement is the union
of open subsets).

\begin{Ex}
A maximal inseparable set need not be connected, for example let
$X$ be the complete graph on the vertex set $V$ with $3< |V|
<\iy$.   The set $V$ is a maximal inseparable subset of $X$.  $V$
also has the property that every pair in $V$ is a cut pair of $X$,
but $V$ is not cyclic.
\end{Ex}
\begin{Lem}\label{L:cyc} Let $S$ be a subset of $X$ with  $|S|>1$.
If every pair of points in $S$ is a cut pair  and  $\{a,b\}$ is a cut pair
separating points of $S$,
then $S \cup \{a,b\}$ is cyclic.
\end{Lem}
\begin{proof}
  It suffices to prove this when $S$ is finite.   Let $c,d \in S$ be
separated
by $\{a,b\}$.  By \fullref{L:cross} , $X- \{c,d\}$ has exactly
two components, $X-\{a,b\}$ has exactly two components, and there
are continua $M_{a,c}$, $M_{a,d}$ $M_{b,c}$, $M_{b,d}$ whose union
is $X$ such that  $M_{a,c}\cap M_{a,d} =\{a\}$, $M_{b,c}\cap
M_{b,d} =\{b\} $, $M_{a,c} \cap M_{b,c} =\{c\}$,  $M_{a,d} \cap
M_{b,d} =\{d\}$.

Now let $e \in S -\{a,b,c,d\}$.    We may assume $e \in M_{a,c}$.  Now
$\{a,b\}$ separates the cut pair $\{d,e\}$ and so by \fullref{L:cross},
$\{d,e\}$ separates $a$ from $b$.  It follows that $e$ is a cut point of the
continuum $M_{a,c}$.  Thus there exist continua $M_{a,e}\ni a$ and
$M_{e,c}\ni c$ such that
$M_{a,e} \cup M_{e,c} = M_{a,c}$ and $M_{a,e} \cap M_{e,c} =\{e\}$.  The set
$\{a,b,c,d,e\}$ is now known to be cyclic.
Continuing this process, we see that $S \cup \{a,b\}$ is cyclic.
\end{proof}
\begin{Cor}
If $S\subset X$ with $|S|>1$  and $S$ has the property that every
pair of points in $S$ is a cut pair, then either $S$ is inseparable or $S$
is cyclic.
\end{Cor}

\begin{Def} By Zorn's Lemma, every cyclic subset of $X$ is contained in a
maximal
cyclic subset.  A maximal cyclic subset with more than two elements is
called a {\em necklace\/}.  In particular, every separable cut pair is
contained in a necklace.
\end{Def}
\begin{Lem} \label{L:twototango} Let $S$ be a cyclic subset of $X$, a
continuum without cut points. If $S$ separates the point $x$ from
$y$ in $X$, then there is a cut pair in  $ S$ separating $x$
from $y$.
\end{Lem}
\begin{proof} Suppose not, then for any finite subset $\{x_1,\dots x_n\}
\subset S$ with $M_1,\dots M_n$ the cyclic decomposition of $X$ by
$\{x_1,\dots x_n\}$, $x$ and $y$ are contained in the same element $M_i$ of
this cyclic decomposition.
There are two cases.

In the first case, we can find a strictly nested intersection of
(cyclic decomposition) continua  $C\ni x,y$,   with the property
that $|C \cap S|\leq 1$.  Nested intersections of continua are
continua, so $C$ connected, and similarly using \fullref{L:cross}, $C- (C \cap S)$ is connected,  so  $S$ doesn't
separate $x$ from $y$.

In the second case there is a cut pair $\{a,b\}\subset S$ and
continua $M,N$ such that $ x,y\in M$, $N \cup M =X$, $N \cap M
=\{a,b\}$ and $ S \subset N$.  It follows that $\{a,b\}$ separates
$x$ from $y$ in $M$, so there exist continua $Y,Z$ with $M=Y\cup
Z$ where $Y\cap Z =\{a,b\}$, $y \in Y$ and $x \in Z$. Since $(N\cup
Y) \cap Z =\{a,b\}$ and $(N\cup Y) \cup Z =X$, it follows that
$\{a,b\}$ separates $x$ from $y$ in $X$.
\end{proof}

\begin{Def} Let $S$ be a necklace of $X$.  We say
$y,z \in X-S$ are $S$ equivalent, denoted $y\sim_S z$, if for every cyclic
decomposition $M_1, \dots M_n$ of $X$ by $\{x_1, \dots x_n \} \subset S$,
both $y,z \in M_i$ for some $1\le i \le n$.  The relation $\sim_S$ is
clearly an equivalence relation on $X-S$.  By \fullref{L:twototango}, if
$y,z$ are separated by $S$ then $y \not \sim_S z$, but the converse is
false.

The closure (in $X$) of a ($\sim_S$)--equivalence class of $X-S$ is called a
{\em gap\/} of $S$.  Notice that every gap is a nested intersection of
continua, and so is a continuum.
Every inseparable cut pair in $S$ defines a
unique gap. The converse is true if $X$ is locally connected, but false in
the nonlocally connected case.

Let $s \in S$.    Choose distinct $x,y \in S-\{s\}$, and take the
cyclic decomposition $M_1,M_2,M_3$ of $X$ by $\{s,x,y\}$ with
$M_1 \cap M_3 =\{s\}$.  For each $i$, take a copy $\hat M_i$ of
$M_i$. Let $\hat M$ be the disjoint union of the $\hat M_i$.  For
$i =1,2,3$ let $s_i,  y_i, x_i$ be the points of $\hat M_i$ which
correspond to $s,x,y$ respectively whenever they exist (for
instance  there is no $s_2$ since $s \not \in M_2$).  Let $\hat X
$ be the quotient space of $\hat M$ under the identification $y_i =y_j$
and $x_i = x_j$ for all $i,j$.   The metrizable continuum $\hat X$
is clearly independent of the choice of $x$ and $y$. The obvious
map $q\co \hat X \to X$ is one to one except that  $\{s_1,s_3\} =
q^{-1}(s)$. We will abuse notation and refer to points of
$X-\{s\}$ as points of $\hat X -\{s_1,s_3\}$ and vice versa.

The cut points of the continuum $\hat X$ are exactly $S-\{s\}$.
Consider the cut point pretree $\cP$ for $\hat X$. By \fullref{L:singing},  the cut points of $\hat X$ will be  exactly the
singleton equivalence classes  in $\cP$ other than $\{s_1\}$ and
$\{s_3\}$. The closures of nonsingleton equivalence classes in
$\cP$ are exactly the  gaps of $S$. {\em Thus  every gap of $S$
has more than one point}. The cut point real tree $T$ is in this
case an arc (see \fullref{L:cross}), so there is a linear order
on $\cP$ corresponding to the pretree structure. Let $A \in \cP$
be a nonsingleton equivalence class (so $\bar A \subset X$ is a
gap of $S$) with $s \not \in \bar A$. Let $U = \{ x \in \hat X:
[x]< A \}$ and let $B = q(\bar U \cap \bar A)$. Similarly let $O=
\{x \in \hat X:A < [x]\}$ and let $C= q(\bar O \cap \bar A)$. The
two closed sets $B$ and $C$ are called the {\em sides\/} of the gap  $\bar
A$.  Notice that $\bd A = C \cup B$. Since $X$ has no cut points $B$ and $C$ 
are nonempty.
\end{Def}
\begin{Def}
Let $D$ be a gap of $S$ with sides $B$ and $C$. If $B\cap C=\emptyset $,
then we say $D$ is  a {\em fat\/} gap of $S$.  Each fat gap  is a continuum 
whose boundary is the disjoint union of its sides.  It follows that every 
fat gap has nonempty interior.  Distinct fat gaps of $S$ will have disjoint 
interiors.
Since the compact metric space $X$ is Lindel\"of (every collection of 
nonempty disjoint open sets is countable), $S$
has only countably many fat gaps.  If $X$ is locally connected then there
are only fat gaps because the sides of a gap form (with local connectivity) an
inseparable cut pair.
\end{Def}
Consider the following example where $X$ is a continuum in  $\R^2$
containing a single necklace $S$ and five gaps of $S$.  The three
solid rectangles are fat gaps, and the two thin gaps  are limit
arcs of Topologist's sine curves.

\begin{center}
\labellist
\pinlabel $X$ at 71 48
\endlabellist
\includegraphics[scale=2.2]{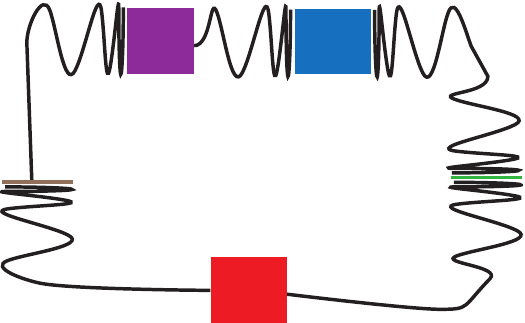}
\end{center}
\begin{Lem} The union of the sides of a gap of $S$ is  a nonsingleton
inseparable set.
\end{Lem}
\begin{proof}
Take $A$, $U$,$O$, $B$,  $C$,  $s$ and $q\co \hat X \to X$  as above.
We show that $B \cup C$ is  a nonsingleton inseparable set.
Suppose that $B\cup C = \{b\}$.  Then $\bd A = \{b\}$ and since gaps are not
singletons, $b$ is a cut point of $X$.  Thus $B \cup C$ is not a singleton.

Now suppose that $d,e \in B\cup C$ and $\{r,t\}$ is a cut pair separating
$d$ and $e$.  Let $E=q(O) \cup q(U)$.  Since $O$ and $U$ are nested unions
of connected sets they are connected and since
$s \in q(O) \cap q(U)$, $E$ is connected.  Thus for any $P\subset X$, with
$E \subseteq P \subseteq \bar E$, $P$ is connected. Since $d,e \in \bar E$
it follows that $\{r,t\}$ must separate $E$.   Since the gap $\bar A\ni d,e$
is connected  it follows that $\{r,t\}$ must separate $\bar A$.  Notice that
$E \cap \bar A$ consists of sides of $A$ which are points of $S$, so $|E\cap
\bar A|\le 2$ and  if $|E\cap \bar A |=2$ then $E \cap \bar A =\{e,d\}$.  It
follows that $\{r,t \} \not \subset E\bar A$.

First consider the case where one of $\{r,t\}$, say $r$ is in $E\cap \bar
A$.
It follows that $r \in S$ is one of the sides of $\bar A$, say $\{r\} =B$.
Thus $d,e \in C$, the other side of $\bar A$.  Since $q(O)$ is connected and
its closure contains $C$, it follows that $t \in q(O) \cap S$. Since $B$ is
not a point of $S$, there are infinitely many  elements $u \in S$ such that
$\{r,u\}$ separates $t$ from $C$. Replacing $t$ with such an $u$, we may
assume that $r$ are $t$ are not inseparable and so $X- \{r,t\}$ has exactly
two components by \fullref{L:cross}.  One of these components will contain
$s$ and the other will contain $\bar A -\{r\}$.  Thus $\{r,t\}$ doesn't
separate $d$ from $e$. Contradiction.

Now we have the case where $\{r,t\} \cap (E\cap \bar A)= \emptyset$.  It
follows that one of them (say $r$)
is a cut point of $E$ and the other $t$ is a cut point of $\smash{\bar A}$.   Since
$r$ is a cut point of $E$, it follows that $r \in S$, and since $r$ is not a
side of $\bar A$, with no loss of generality $r=s$.  Thus $t$ is a cut point
in $\hat X$ and so $t \in S$.  But $S \subset E$ and $t \not \in E$.
Contradiction.
\end{proof}
\begin{Cor} \label{C:character}
Let $X$ be a continuum without cut points.
Suppose that for every pair of points $c,d\in X$ there is a pair
of points $a,b$ that separates $c,d$. Then $X$ is homeomorphic to
the circle.
\end{Cor}
\begin{proof}
Let $S$ be a necklace of $X$. Using \fullref{L:cross}, we can show that $S$ is infinite and in fact any
two points of $S$ are separated by a cut pair in $S$.
If $X-S \neq \emptyset $ then there is a gap $A$ of $S$. The union of sides
of the gap $A$ is a nonsingleton inseparable subset of $X$.  There are no
nonsingleton inseparable subsets of $X$, so $S=X$.
Thus  $X$ is homeomorphic to the circle.  This follows from 
\cite[Theorem 2-28, page 55]{HOC-YOU} (our  \fullref{L:maptocircle} also proves this).
\end{proof}

\begin{Thm} \label{L:maptocircle}
Let $S$ be a necklace of $X$. There exists a continuous surjective
function $f\co X \to S^1$, with the following properties:
\begin{enumerate}
\item The function $f$ is one to one on $S$. \item The image of a
fat gap of $S$ is a nondegenerate arc of $S^1$. \item For $x,y\in
X$  and  $a,b\in S$:
\vspace{10pt}
\begin{enumerate}
\item
If  $\{f(a), f(b)\}$ separates $`f(x)$ from $`f(y)$  then $\{a,b\}$ separates
$x$ from $y$.
\item  If $x \in S$ and $\{a,b\}\subset S$ separates $x$ from $y$, then
$\{f(a),
f(b)\}$ separates $f(x)$ from $f(y)$.
\end{enumerate}
\end{enumerate}
The function $f$ is unique up to homotopy and reflection in $S^1$.
In addition, if $G$ is a group acting by homeomorphisms  on $X$
which stabilizes $S$,  then the action of $G$ on $S$ extends to an
action of $G$ on $S^1$.
\end{Thm}

\begin{proof}
We use the strong  Urysohn Lemma \cite[4.4 Exercise 5]{MUN} If $A$
and $B$ are disjoint closed $G_\delta$ subsets of a normal space $Y$, then
there is a continuous $f\co Y \to [0,1]$ such that $f^{-1}(0) = A$ and
$f^{-1}(1) = B$.   In a metric space, all closed sets are $\smash{G_\delta}$.
Since $X$ has a countable basis, the subspace $S$ has a
countable dense subset
$\smash{\hat S}$.  Since the fat gaps of $S$ are countable, the collection $R$ of all 
sides of fat gaps of $S$ is countable. Let $\{s_n:n\in \N \}= \hat S \cup 
R$.  Notice that now some of the elements of $\{s_n: n \in \N\}$ are points 
(singleton sets) of
$S$ and some of them are sides  of gaps (and therefore closed sets of $X$). 
In particular all inseparable cut pairs of $S$ are in $\{s_n:n \in \N\}$.

For the remainder of this proof, we will maintain the useful fiction that
each element of $\{s_n:n \in \N\}$ is a point (which would be true if $X$ 
were locally
connected), and leave it to the reader (with some hints) to check the
details for the nonsingleton sides of gaps.

  Notice that the elements of $\{s_n: n \in \N\}$ are pairwise disjoint.

We inductively construct the map $f$.  We take as $S^1$, quotient space of
the interval $[0,1]/(0=1)$ with $0$ identified with $1$.
Since $\{s_1,s_2,s_3\}$ is cyclic, there exist cyclic decomposition
$M_1,M_2,M_3$ of $X$ with respect to $\{s_1,s_2,s_3\}$.  We define the map $f_3 \co M_
1\to [0,\unfrac 13]$ by $f_3(s_1) =0$, $f_3(s_2) = \unfrac 1 3$ and then extend
to $M_1$ using the strong Urysohn Lemma so that $f_3^{-1}(0)
=\{ s_1\}$ and
$ f_3^{-1}(\unfrac 1 3) =\{s_2\}$.  Similarly we define the continuous map $
f_3\co M_2 \to [\unfrac 1 3, \unfrac 2 3] $ such that $ \smash{f_3^{-1}}( \unfrac 1 3)
=\{s_2\}$ and $ \smash{f_3^{-1}}(\unfrac 2 3) = \{s_3\}$.  Lastly we define $
f_3\co M_3 \to [\unfrac 2 3, 1] $ such that $ \smash{f_3^{-1}}( \unfrac 2 3) =\{s_3\}$
and $ \smash{f_3^{-1}}(1) = \{s_1\}$.  Since $0=1$ we paste to get the
function $f_3\co X \to S^1$.

Now inductively suppose that we have $N_1, \dots N_k$ a cyclic decomposition
of $X$ with respect to
$\{s_i : i \le k\}$ (when the $s_i$ are sides of gaps the definition of
cyclic decomposition will be similar), and a map $f_k\co X \to S^1$ such that
$f_k(N_j) = [f_k(s_p),f_k(s_q)]$ for each $1\le j
\le k$,   where $\bd N_j =\{s_p,s_q\}$ and
$q,p \le k$, satisfying $f_k^{-1}(f(s_j))=\{s_j\}$ for all $j\le k$.
If $s_{k+1} \in N_j$ with $\bd N_j =\{s_p,s_q\} $ then there exists continua
$A,B$ such that $A\cup B = N_j$, $A \cap B = \{s_{k+1}\}$, $s_p \in A$ and
$s_q \in B$ (in the case where $s_{k+1}$ is the side of a gap, then  one of
$A$, $B$ will
be a nested union of continua, and the other will be a nested intersection).
Using the strong Urysohn Lemma, we
define $f_{k+1}\co N_j \to [f_k(s_p),f_k(s_q)]$ such that
$\smash{f_{k+1}^{-1}}(f_k(s_p))=\{s_p\}$, $\smash{f_{k+1}^{-1}}(f_k(s_q))=\{s_q\}$,
$f_{k+1}^{-1}(
\upnfrac {s_q +s_p}2)=\{s_{k+1}\}$,  $f_{k+1}(A) =[f_{k+1}(s_p),
f_{k+1}(s_{k+1})]$
and lastly $f_{k+1}(B) =[f_{k+1}(s_{k+1}), f_{k+1}(s_{q})]$.  We define $f_{k+1}$
to be equal to $f_k$ on $X- N_j$ and by pasting we obtain $f_{k+1}\co X \to
S^1$.
By construction, the sequence of functions $f_k$ converges uniformly to a
continuous function
  $f\co X \to S^1$.   Property (2) follows from the construction of $f$.

  For uniqueness, consider $h\co X \to S^1$ satisfying these properties.  Since
the cyclic ordering on $S$ implies that $f\co S\to S^1$ is unique up
to isotopy and reflection \cite{BOW1}, we may assume that $h$ and
$f$ agree on $S$. Thus for any fat gap $O$ of $S$, we have $h(O) =
f(O)=J$, an interval.  Since  $f$ and $h$ agree on the sides of
$O$, which are sent to the endpoints of $J$, we simply straight
line homotope $h$ to $f$ on each fat gap.  Clearly after the
homotopy they are the same.

  The action of $G$ on $S$ gives an action on $\overline{f(S)} \subset S^1$,
which preserves the cyclic order.  Thus by extending linearly on the
complementary intervals, we get an action of $G$ on $S^1$.  This action has
the property that for any $g \in G$, $f\circ g \simeq g \circ f$.
\end{proof}

\begin{No}Let $X$ be a continuum without cut points.  We define $\cR\subset
2^X$ to be the collection of all necklaces of $X$, all maximal
inseparable subsets of $X$, and all inseparable cut pairs of $X$.
For the remainder of this section, $X$ is fixed.
\end{No}
\begin{Lem}\label{L:cot} Let $E$ be a nonsingleton subcontinuum of $X$.
There exists $Q \in \cR$ with $Q \cap E \neq \emptyset$.
\end{Lem}
\begin{proof}
Let $c,d\in E$ distinct.  If $\{c,d\}$ is an inseparable set, then there is
a maximal inseparable set $D \in \cR$ with $c,d \in D$.

If not then there is a cut pair $\{a,b\}$ separating $c$ from $d$.  It
follows that $E\cap \{a,b\}\neq \emptyset$.  Either $\{a,b\}$ is inseparable or there is a necklace $N \in \cR$
with $\{a,b \} \subset N$, and so $E \cap N \neq \emptyset$.
\end{proof}
\begin{Thm}\label{T:cint} Let $X$ be a continuum without cut points. If $S,
T\in \cR$ are distinct then $|S\cap T|<3$ and if $|S\cap T| =2$, then $S
\cap T$ is an inseparable cut pair.
\end{Thm}
\begin{proof}
If $S$ or $T$ is an inseparable cut pair, then the result is trivial. We are
left with three cases.

First consider the case where $S$ and $T$ are necklaces of $X$
Suppose there are distinct $a,b,c \in S \cap T$.  Since $S,T$ are distinct
necklaces, there exists $d\in S-T$.    Since $\{a,b,c,d\}
\subset S$ is cyclic, renaming $a,b,c$ if needed, we have  $X= A\cup B \cup
C \cup D$ where $A,B,C,D$ are continua and $A\cap B= \{b\}$, $B\cap C =
\{c\}$, $C\cap D =\{d\}$, and $D\cap A = \{a\}$, and all other pairwise
intersections are empty.    Thus $\{b,d\}$ separates $a$ and $c$, points of
$T$.  It follows by
\fullref{L:cyc} that $d \in T$. This contradicts the choice of $d$ so $|S
\cap T|<3$.
Now suppose we have distinct $a,b \in S \cap T$. If $\{y,z\}$ is a cut pair
separating $a$ from $b$ in $X$ then, by \fullref{L:cyc}, $\{y,z\} \subset
T$ and $\{y,z \} \subset S$, so
$|S\cap T| >3$.  This is a contradiction, so $\{a,b\}$ is an inseparable cut
pair.

Now consider the case where $S$ and $T$ are maximal inseparable
subsets of $X$.  Since $S$ and $T$ are distinct maximal
inseparable sets, there exist $y \in S $, $z \in T$ and a cut pair
$\{a,b\}$ separating $y$ from $z$.  It follows that $y \not \in S$
and $z \not \in T$. Thus $X = C \cup D$ where $C$ and $D$ are
continua,  $y \in C$, $z \in D$ and $C \cap D = \{a,b\}$. By
inseparability, $S \subset C$ and $T \subset D$. Clearly $S\cap T
\subset C \cap D =\{a,b\}$. If $S\cap T = \{a,b\}$ then $\{a,b\}$
is inseparable.

Finally consider the case where $S$ is a necklace of $X$, and
$T$ is a maximal inseparable set of $X$.  By definition, every cyclic subset
with more than three elements is not inseparable.  It follows that $|S\cap
T|
<4$.   The only way that $|S \cap T| =3$ is if $S = T$ which is not allowed.
If $|S\cap T|=2$ then $S\cap T$ is inseparable (since $T$ is) and
cyclic (since $S$ is) and therefore $S\cap T$ is an inseparable cut pair.
\end{proof}

\begin{Lem} \label{L:nosep} If $S,T \in \cR$, then $S$ doesn't separate
points of $T$.
\end{Lem}
\begin{proof}
Suppose that $r,t \in T-S$ with $S$ separating $r$ and $t$. First
suppose that $S$ is cyclic (so $S$ is a necklace or an inseparable
cut pair).  In this case by \fullref{L:twototango}, there exists
a cut pair $\{a,b\} \subset S$ such that $\{a,b\}$ separates $r$
from $t$.

If $\{r,t\}$ is a cut pair, then by  \fullref{L:cross}, $a$ and
$b$ are separated by $\{r,t\}$, so $S$ is not  an inseparable
pair. Thus $S$ is a necklace and it follows by \fullref{L:cyc}
that $r,t \in S$ (contradiction). Thus $\{r,t\}$ is not a cut
pair, and so $T$ is a maximal inseparable set, but $\{a,b\}$
separates points of $T$ which is a contradiction.

We are left with the case where $S$ is a maximal inseparable set.

If $T$ is also a maximal inseparable set, then there is a cut pair $A$
separating a point of $S$ from a point of $T$.  Thus there exist  continua
$N$ and $M$ such that $N\cup M =X$, $N \cap M =A$ and, since $S$ and $T$ are
inseparable, we may assume  $T\subset N$ and $S \subset M$.  Since $A$
doesn't separate points of $T$,  and $S \cap N \subset A$, it follows that
$S$ doesn't separate points of $T$.

Lastly we have the case where $S$ is maximal inseparable, and $T$ is cyclic.
Thus $\{r,t\}$ is a cut pair.
So there exist continua $N,M$ such that $N\cup M =X$ and $N \cap M
=\{r,t\}$.   Since $S$ is maximal inseparable, $S$ is contained in one of
$N$ or $M$ (say $S \subset M$).  However $r,t \subset N$ and since
$r,t\not \in S$, we have $S \cap N = \emptyset $.  Thus $S$ doesn't separate $r$
from $t$.
\end{proof}

\begin{Def}
We now define a symmetric betweenness relation on $\cR$  under which $\cR$ is
a pretree.
Let $R,S,T $ be distinct elements of $\cR$.  We say $S$ is between $R$ and
$T$, denoted $RST$ or $TSR$, provided:
\begin{enumerate}
\item \label{cutpair}   $S$ is an inseparable cut pair and $S$ separates a
point of $R$ from a point of $T$.
\item   $S$ is not an inseparable pair and:
\vspace{10pt}
\begin{enumerate} \item
$R \subset S$, so $R$ is an inseparable cut pair,   and $R$ isn't between
$S$ and $T$ (see case \eqref{cutpair}).
\item \label{partof} $S$ separates a point of $R$ from a point of $T$, and
there is no  cut pair $Q \in \cR$ with $RQS$ and $TQS$ (see case
\eqref{cutpair}).
\end{enumerate}
\end{enumerate}
For $R,T \in \cR$ we define the open interval $(R,T) = \{S\in
\cR:\, RST\}$. We  now defined the closed interval $[R,S]= (R,S)
\cup \{R,T\}$ and we define the half-open intervals analogously.
We will show that $\cR$ with this betweenness relation forms a
pretree \cite{BOW5}. Clearly for any $R,S \in \cR$,  by
definition $[R,S] = [S,R]$ and $R \not \in (R,S)$.
\end{Def}

Consider the following example in  
where $X \subset \R^2$ is the union of 6
nonconvex quadrilaterals (meeting  only at vertices) and a Topologist's
sine curve limiting up to one of them.  
There are two
necklaces, one being the Topologist's sine curve and the other consisting of
the four green points.
The cut pair tree $T$ is obtained from $\cR$ by connecting the dots
(definition to be given later).
\begin{center}
\labellist
\pinlabel $X$ at 200 39
\endlabellist
\includegraphics[scale=1.2]{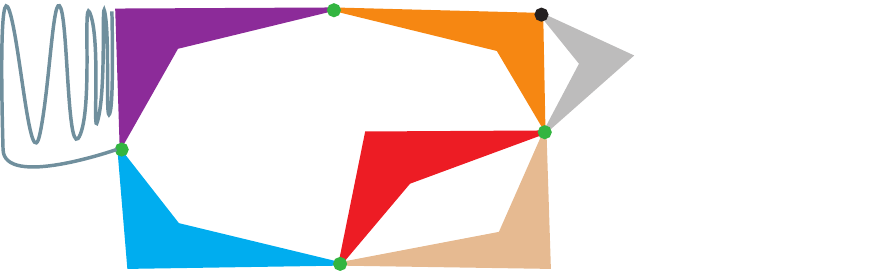}
\labellist
\pinlabel {$\text{The pretree } \cR$} at 195 35
\endlabellist
  \includegraphics[scale=1.2]{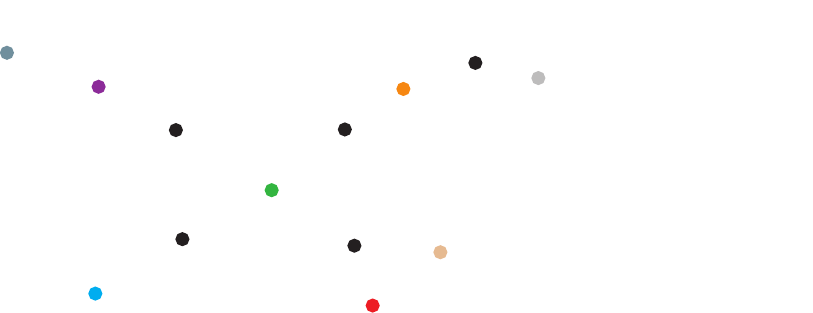}
\labellist
\pinlabel {$\text{The cut pair tree } T$} at 195 35
\endlabellist
   \includegraphics[scale=1.2]{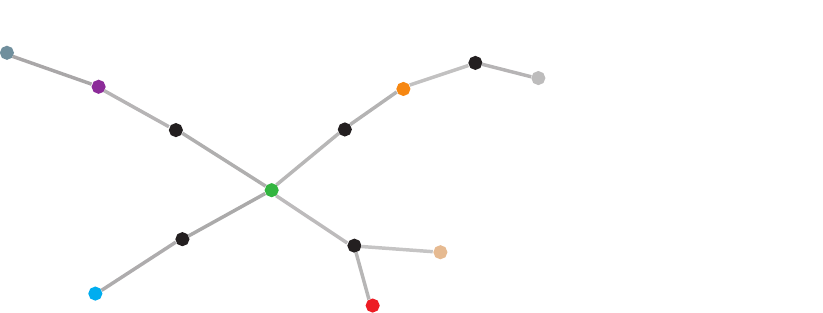}
\end{center}

\begin{Lem}  For any $R,S,T \in \cR$, we have that $[R,T] \subset [R,S]\cup 
[S,T]$.
\end{Lem}
\begin{proof}
We may assume $R,S,T$ are distinct.
Let $ Q \in (R,T) $ with $Q \neq S$.

   If $Q$ is an inseparable pair  then
$Q$ separates a point $r\in R$ from a point $t\in T$.
Thus there exist continua $N$,$M$ such that $N\cup M =X$, $N\cap M = Q$,
$r\in N$ and $t\in M$.   Since $S \not \subset Q$, either $(S-Q) \cap N \neq
\emptyset$ implying $Q \in (S,T)$, or
$(S-Q) \cap M \neq \emptyset $ implying $Q \in (R,S)$.

Now consider the case where $Q$ is not an inseparable pair.

Suppose that one of $R$,$T$  (say $R$) is contained in $Q$,  so $R\subset Q$
is an
inseparable cut pair and $R\not \in (Q, T)$.  If $R \not \in (Q, S)$ then by
definition $Q \in (R,S)$ as required.
If on the other hand $R \in (Q,S)$  then there exist continua $N, M$ and
$q\in Q-R$ and $s \in S-R$ such that $q\in N$, $s \in M$, $N\cup M =X$ and
$N \cap M = R$.  Since $R\not \in (Q, T)$, it follows that  $(T-R) \subset
N$.   If $T \subset Q$, then since $T\neq R$, it follows that $T \not \in
(Q,S)$, and so $Q \in (T,S)$.  If $T\not \subset Q$, then there is $t \in
(T-Q)\subset  (T-R) \subset N$ and it follows that $Q$ separates $s$ from
$t$ since $R$ separated them, thus $Q \in (S,T)$ as required.

We are now left with the case where $Q$ is not an inseparable pair, $R
\not \subset Q$, and $T \not \subset Q$ (see Property (2)\eqref{partof} of the betweenness relation).  Thus by
definition there exists $r \in R-Q$, $t\in T-Q$ and disjoint continua $M,N$
with $r \in M $ and $t \in N$, $N\cup M = X$ and $N \cap M \subset Q$.  If
$S \not \subset Q$ then there exists $s\in S-Q$ and either $s \in M$ in
which case $Q \in (S,T)$ or $s \in N$ in which case $Q \in (S,R)$.
If on the other hand $S \subset Q$, then by Property (2)\eqref{partof}, either $S \not
\in (Q,R)$ implying $Q \in (S,R)$ or $S \not \in (Q,T)$ implying $Q \in
(S,T)$.
\end{proof}
\begin{Lem} For any $R,T \in \cR$, if $S \in (R,T)$ then $R \not \in
(S,T)$.
\end{Lem}
\begin{proof}
First consider the case where $S$ is an inseparable cut pair.  We have $r
\in R-S$, $t \in T-S$ and continua $N\ni r$ and $M\ni t $ such that $N\cup M
=X$, $N \cap M = S$.  In fact by \fullref{L:nosep}
$R \subset N$ and $T \subset M$.

  If $S \not \subset R$, then $|R\cap S|<2$. Since $X$ has no cut points, no
point in $S$ is a cut point of $M$, so $M-R$ is connected. Thus  $R$ doesn't
separate
$S$ from $T$, so $R \not \in (S,T)$.

If $S \subset R$,  then by definition since $S \in (R,T)$ then $R \not \in
(S,T)$.

Now consider the case where $S$ is not an inseparable pair.
If  $R \subset S$, then  $R$ is an inseparable pair and $R \not \in (S,T)$
as required.
We may now assume that $R \not \subset S$.
If $T \subset S$, then $T$ is an inseparable pair and $T\not \in (R,S)$.  By
\fullref{L:nosep}
$R$ cannot separate a point of $T$ from a point of $S$, since $R\not \subset
S$, it follows that $R \not \in (S,T)$.

We are left with case (2,b) , so $S$ separates a point $r \in R-S$
from a point $t \in T-S$.  Thus there exists continua $M,N$ with $r \in M$,
$t\in N$, $N \cup M =X$ and $N \cap M \subset S$. In fact by \fullref{L:nosep}
$R \subset M$ and $T \subset N$.
Since $X$ has no cut points and $|R \cap S|<2$, then $N-R$ is connected, and
so $R \not \in (S,T)$.
\end{proof}

\begin{Def} We say distinct $R,S \in \cR$ are adjacent if $(R,S) =
\emptyset$.
\end{Def}
\begin{Lem}
If $R,S \in \cR$ are adjacent then $R \subset S$, $S \subset R$,  or
(interchanging if need be) $R$ is a necklace and  $S$ is maximal inseparable
with $[\bar R - R]\cap S \neq \emptyset$.
\end{Lem}
\begin{proof}
We need only consider the case where $R,S$ are adjacent and neither is a
subset of the other.

First consider the case where  one of $R$, $S$ (say $S$) is an inseparable
set.  There is no maximal inseparable set containing both $R$ and $S$, so
there exists $r\in R-S$ and  cut pair $A$ separating $r$ from a
point of $S$.  Notice that  $A$ is contained in some necklace  $T$.
Since $A,T \not \in  (R,S)$, it follows that  $T=R$ and that $S$ is maximal
inseparable.

  Let $G$ be the gap of $R$ with $S\subset G$.  Let $Q$ be a side of $G$ and
$p \in Q$.   If $p \not \in S$, then there exists a cut pair $B$ separating
$p$ from $S$.
Since $(R,S)= \emptyset$, $B$ doesn't separate $R$ from $S$.  It follows by
definition of side, that $B$ separates points of $R$ which implies that $B
\subset R$.  This contradicts the fact that $Q$ is a side of the gap
$G\supset S$.   If both sides of $G$ are points, then they form an
inseparable cut pair in $(R,S)$.  Thus they are not both points so $[\bar
R - R]\cap S \neq \emptyset$.

We are left with the case where $R$ and $S$ are each necklaces with more
than $2$ elements.
Again let $G$ be the gap of $R$ with $S \subset G$, and let $Q,P$ be sides
of $G$.  Since $Q \cup P$ is inseparable, there is a maximal inseparable set
$A\supset Q \cup P$.  It follows that
$A \in (R,S)$ which is a contradiction.
\end{proof}
Using the pretree structure on $\cR$, we can put a linear order
(two actually) on any interval of $\cR$. We recall that the order
topology on a linearly ordered set $I$ is the topology with
basis $I_y=\{x:x>y \}, J_y=\{x:x<y \}$ and $K_{y,z}=\{x:z<x<y
\}$ where $y,z$ range over elements of $I$.
The suspension of a Cantor set
is a continuum with uncountably many maximal inseparable
sets, but this doesn't happen for inseparable cut pair and necklaces.
\begin{Lem} \label{L:count} Only countably many elements of $\cR$  are 
inseparable pairs or necklaces.
\end{Lem}
\begin{proof}
We first show that any interval $I$ of $\cR$ contains only countably many 
necklaces and inseparable pairs.

Let $Q$ be the set of all cut pairs in $I$  which have more than 
two complementary components, union the set of necklaces in $I$.
Let $A \in Q$.
\begin{itemize}
\item
If $A$ is a cut pair, then since $X-A$ has more than two components, and 
$\smash{\bigcup Q}$ will intersect two of the components, $A$ will separate  $(\bigcup 
Q)-A $ from some other point of $X$.  Using Lemma 3, we find subcontinua 
$Y,Z$ of $X$ such that  $ Y \cup Z=X$, $Y \cap Z = A$ where $(\bigcup Q)-A 
\subset Y$. We define the open set $U_A= Z-A$
\item If $A$ is a necklace, then $|A|>2$ and there is a cut pair $\{a,b \} 
\subset A$ which doesn't separate $(\bigcup Q)-A$.  Using Lemma 3, we find 
subcontinua $Y,Z$ of $X$ such that  $ Y \cup Z=X$, $Y \cap Z = A$ where 
$(\bigcup Q)- A\subset Y$. We define the open set $U_A= Z-\{a,b\}$
\end{itemize}
Notice that for any $A,B \in Q$, $U_A \cap U_B= \emptyset$. Since $X$ is 
Lindel\"of, the collection $\{U_A: A \in Q\}$ is countable and therefore $Q$ 
is countable.

It is more involved to show that inseparable cut pairs $\{a,b\}$ of $I$
such that $X-\{a,b\}$ has 2 components are countable.
Let $S$ be the set of inseparable cut pairs $\{a,b\}$ in
$I$ such that $X-\{a,b \}$ has 2 components. We argue by
contradiction, so we assume that $S$ is uncountable.

Let $\{a,b\}$ be a cut pair of $S$ and let $C_L,C_R$ be the
components of $X-\{a,b\} $.

We say that $\{a,b\}$ is a limit pair if there are inseparable cut
pairs $\{a_i,b_i\}$ and $\{a_i',b_i'\}$ in $S$ such that
$\{a_i,b_i\}\subset C_L$, $\{a_i',b_i'\}\subset C_R$, and for each
limit pair $\{ a',b' \}\ne \{a,b \}$ of $S$ one of the two
components of $X-\{a',b'\} $ contains at most finitely many
elements of the sequences $\{a_i,b_i\}$ and $\{a_i',b_i'\}$.

  We claim that there are at most
countable pairs in $S$ which are not limit pairs. Indeed if $\{
a,b \}$ is not a limit pair and $I=[x,y]$ let $C_L,C_R$ be the
components of $X-\{ a,b \}$ containing, respectively, $x,y$ ($L,R$
stand for left, right). Since $\{ a,b \}$ is not a limit pair for
some $\epsilon
>0$ one of the four sets
$$C_L\cap B_{\epsilon }(a),\ C_R\cap B_{\epsilon
}(a),\ C_L\cap B_{\epsilon }(b,)\text{ or }C_R\cap B_{\epsilon }(b)$$
intersects the union of all cut pairs of $S$ at either $a$ or
$b$.

We remark now that for fixed $\epsilon >0$ there are at most
finitely many pairs $\{a,b \}$ in $S$ such that (say) $C_L\cap
B_{\epsilon }(a)$ intersects the union of all cut pairs of $S$
in a subset of $\{a,b \}$. Indeed if we take all pairs $\{ a,b \}
$ with this property the balls $B_{{\punfrac\epsilon 2} }(a)$ are
mutually disjoint so there are finitely many such pairs. The same
argument applies to each one of the three other sets $C_R\cap
B_{\epsilon }(a),\ C_L\cap B_{\epsilon }(b,)$ or $C_R\cap
B_{\epsilon }(b)$. This implies that non limit cut pairs are
countable.

So we may assume $S$ has uncountably many limit pairs. Let
$\{c,d \}$ be a limit pair in $S$, let $C_L,C_R$ be the
components of $X-\{c,d \}$ and let $\{c_i,d_i \}\in \bar C_L$,
$\{c_i',d_i' \}\in \bar C_R$ be sequences of distinct pairs in
$S$ provided by the definition of limit pair.

Let $\smash{C_R^i}$ be the component of $X-\{c_i,d_i\}$ containing $c,d$.
We claim that there is an $\epsilon $ such that for all $i$ there
is some $x_i\in C_L\cap C_R^i$ with $d(x_i,\{ c,d \})> \epsilon $.
Indeed this is clear if the accumulation points of the sequences
$c_i$ and $d_i$ are not contained in the set $\{c,d \}$. Otherwise
by passing to a subsequence and relabelling, if necessary, we may
assume that either $c_i\to c,d_i\to d$ or both $c_i,d_i$ converge
to, say, $c$.

In the first case we remark that there is a component $C_i$ of
$X-\{c,d,c_i,d_i\}$, such that its closure contains both $c,d_i$
or both $d,c_i$. Indeed otherwise $\{c,d,c_i,d_i\}$ is a cyclic
subset which is impossible since we assume that $\{ c_{j},d_{j}
\}\ (j>i)$ are all inseparable cut pairs.

Since $d_i\to d$ and $c_i\to c$ there exists $\epsilon >0$ and
$x_i\in C_i$ such that $d(x_i,c)>\epsilon$ and $d(x_i,d)>\epsilon $
for all $i$.

In the second case we remark that since $c$ is not a cut point
there is some $e >0$ such that for each $i$ there is a component
$C_i$ of $C_L\cap \smash{C_R^i}$ with diameter bigger than $e $. It
follows that there is an $\epsilon >0$ and $x_i\in C_i$ so that
$d(x_i,c)>\epsilon ,d(x_i,d)>\epsilon $ for all $i$.

 By passing to subsequence we may assume that $x_i$
converges to some $x_L\in C_L$. Clearly $d(x_L,c)\geq \epsilon
,d(x_L,d)\geq \epsilon $. It follows that $d(x_L,C_R)>0$.

We associate in this way to a limit pair $\{c,d
\}$ in $S$ a point $x_L$ and a $\delta >0$ such that:
\begin{enumerate}
\item $x_L\in C_L$
\item$d(x_L,C_R) >\delta $
\end{enumerate}

Since there are uncountably many limit pairs in $S$ there are
infinitely many such pairs for which $1,2$ above hold for some
fixed $\delta >0$. But then the corresponding $x_L$'s are at
distance greater than $\delta $ (by property 2 above). This is
impossible since $X$ is compact. Thus $S$ is countable.

Thus for any interval $I$ of $\cR$, the set of necklaces and inseparable cut 
pairs of $I$
is countable.

  Let $E$ be a countable dense subset of $X$.
For any $A$, a necklace with more than one gap or an inseparable pair, there 
exist
$a,b \in E$ separated by $A$.  Thus the interval $I=[[a],[b]]$  contains 
$A$.  There are countably many such intervals, so the set of necklaces with more 
than one gap is countable, and the set of inseparable cut pairs of $X$ is 
countable.

 If a necklace $N$ has less than two gaps, there is an open set $U \subset 
N$. By Lindel\"of, there are at most countably many such necklaces, and thus 
there are at most countably many inseparable cut pairs and necklaces in $X$.
\end{proof}

\begin{Lem} The pretree $\cR$ is preseparable and complete.
\end{Lem}
\begin{proof} Let $[R,W]$ be a closed interval of $\cR$.

We first show that any  bounded strictly increasing sequence in $[R,W]$
converges.
Let  $(S_n)\subset [R,W]$ be strictly
increasing.  Let $C_n$ be the component of $X-S_n$ which contains $R$.  Let
$C= \overline{\bigcup C_n}$.   Clearly $C$ is contained in the closure $Q$ of
the component of $X-W$ which contains $R$, and so $\bd C \subset Q$.
Clearly $\bd C$ is not a point (since it would by definition be a cut point
separating $R$ from $W$).  The set $\bd C$ is inseparable, and so $\bd C
\subset A$ is a maximal inseparable set. It follows that $A \in [R,W]$.
If $S_n \not \to A$, then there is $B \in [R,A)$ with $S_n <B$ for all $n$.
As before, we have $C$ contained in the closure $D$ of the component of
$X-B$ containing $R$.
This would imply that $A \in [R,B]$, a contradiction.
Thus every strictly increasing sequence in $[R,W]$ converges.

We now show that there are only countably many adjacent pairs in
$[R,W]$. We remark that if $A,B$ is an adjacent pair in $\cR $ at
most one of the sets $A,B$ is a maximal inseparable set.
By \fullref{L:count} there are only countably many inseparable cut pairs 
and necklaces in
$X$.  It follows that there are only countably many inseparable pairs in 
$[R,W]$.
\end{proof}

We have shown thus that $\cR $ is a complete preseparable pretree.
By gluing in intervals to adjacent pairs of $\cR $ we obtain a
real tree $T$ as in \fullref {T:countable}.

\begin{Cor} There is a metric on $T$, which preserves the pretree  structure
of $T$, such that $T$ is an  $\mathbb{R} $--tree. The topology so
defined on $T$ is canonical.
\end{Cor}
\begin{proof}
We metrize $T$ as in \fullref {T:Rtree}. We metrize first the
subtree spanned by the set of inseparable cut pairs and necklaces
(which is countable)
and then we glue intervals for the inseparable subsets of $\cR $
which are not contained in this subtree.
\end{proof}
We call this $\mathbb{R} $--tree the \textit{JSJ-
tree\/} of the continuum $X$ since in the case $X=\partial G$ with
$G$ one-ended hyperbolic our construction produces a simplicial
tree corresponding to the JSJ decomposition of $G$.

\section{Combining the two trees}
When $X$ is locally connected, one can combine the constructions of the
previous 2 sections to obtain a tree for both the cut points and the cut 
pairs of a
continuum $X$.  The obvious application would be to relatively hyperbolic 
groups, and we should note that in that setting, the action of the tree may 
be nesting.
 We explain briefly how to construct this tree.

Let $X$ be a Peano continuum and let  $\cP $ be the cut point pretree.

\begin{Lem} Let $A\in \cP$ be a nonsingleton equivalence class of $\cP$.
Then the closure
$\bar A$ is a Peano continuum without cut points.
\end{Lem}
\begin{proof}
We first show that $\bar A$ is a Peano continuum.  Since $\bar A$ is
compact, and $X$ is (locally) arcwise connected, it suffices to show that
$A$ is convex in the sense that every arc joining points of $\bar A$ is
contained in $\bar A$.

Let $a,b$ be distinct points of $ \bar A$ and let $I$ be an arc from $a$ to
$b$. Suppose $d \in I-A$.
Thus either $c$ is a cut point adjacent to $A$, or there is a cut point $c
\in A$ separating $d$ from $A$,
but then $I$ cannot be an arc since it must run through $d$ twice.

Let $a,b,e\in \bar A$.  Since $e$ doesn't separate $a$ from $b$ in $X$,
there is an arc in $X$ from $a$ to $b$ missing $e$.  By convexity, this arc
is contained in $\bar A$, and so $e$ doesn't separate $a$ from $b$ in $\bar
A$.  It follows that the continuum $\bar A$ has no cut points.
\end{proof}
Let $A$ be a nonsingleton equivalence class of  the cut point pretree
$\cP$, and let $T_A$ be the ends compactification (well,  it will not be
compact, but we glue the ends to the tree anyway) of the cut pair tree for
$\smash{\bar A}$.
Since $X$ is locally connected, for any interval $(B,D)\ni A$ there are cut
points $a_1, a_2\in \smash{\bar A}$ with $a_1,a_2 \in (B,D)$.
Not every point of $\smash{\bar A}$ is contained in one of the defining sets of the
cut pair pretree for $\smash{\bar A}$.
Some of the points of $\smash{\bar A}$ are  not contained in a cut pair, or in a
maximal inseparable set with more than two elements, and these appear as
ends of the cut pair tree $\cR_{\bar A}$ for $\smash{\bar A}$.

For each nonsingleton class $A$ of a the cut point tree $T$ we replace $A$
by $T-A$.  The end of the component of $T-A$ corresponding to a cut point
$a_1 \in \smash{\bar A}$ is glued to the minimal point  or end of $T_A$ containing
$a_1$.

To see that this construction yields a tree,  we use the following Lemma.

\begin{Lem} The set of classes of $\cP$ with nontrivial relative
JSJ-tree in any interval of $\cP $ is countable.
\end{Lem}
\begin{proof}
Let $[u,v]$ be an interval of $\cP $ and let $A$ be a class of
$[u,v] $ with nontrivial JSJ-tree. Since $X$ is locally connected,
$\bar A$ contains some cut point of $[u,v]$. If $c$ is a cut point of
$[u,v]$ in $\bar A$ we have that
$c,A$ are adjacent elements of $[u,v]$. But we have shown in
\fullref {T:countable} that there are at most countable such
pairs.
\end{proof}
Clearly any group of homeomorphisms of $X$ acts on this combined tree.

%
%
%
%
%

\section{Group actions}

The $\Bbb R$--trees we construct in the previous sections usually
come from group boundaries and the group action on them is induced
from the action on the boundary, so it's an action by
homeomorphisms. In this section we examine such actions and
generalize some results from the more familiar setting of
isometric actions.

We recall that the action of a group $G$ on an $\Bbb R$--tree $T$
is called non-nesting if there is no interval $[a,b]$ in $T$ and
$g\in G$ such that $g([a,b])$ is properly contained in $[a,b]$. An
element $g\in G$ is called elliptic if $gx=x$ for some $x\in T$.
If $g$ is elliptic we denote by $\fix(g)$ the fixed set of $g$. An
element which is not elliptic is called hyperbolic.
\begin{Lem}
Let $G$ be a group acting on an $\Bbb R$--tree $T$ by
homeomorphisms. Suppose that the action is non-nesting. Then if
$g$ is elliptic $\fix(g)$ is connected. If $g$ is hyperbolic then
$g$ has an ``axis'', ie there is a subtree $L$ invariant by $g$
which is homeomorphic to $\Bbb R$.
\end{Lem}
\begin{proof}
Let $g$ be elliptic. We argue by contradiction. If $A,B$ are
distinct connected components of $\fix(g)$ let $[a,b]$ be an
interval joining them ($a\in A, b\in B$). Then $g([a,b])=[a,b]$.
Since $[a,b]$ is not fixed pointwise there is a $c\in [a,b]$ such
that $g(c)\ne c$. So $g(c)\in [a,c)$ or $g(c)\in (c,b]$. In the
first case $g([a,c])\subset [a,c)$ and in the second
$g([c,b])\subset (c,b]$. This is a contradiction since the action
is non-nesting.

Let $g$ be hyperbolic. If $a\in T$ consider the interval
$[a,g(a)]$. The set of all $x\in [a,g(a)]$ such that $g(x)\in
[a,g(a)]$ is a closed set. If $c$ is the supremum of this
set then there is no $x\in [c,g(c)]$ such that $g(x)\in [c,g(c)]$.
We take $L$ to be the union of all $g^n([c,g(c)])$ ($n\in \Bbb
Z$). Clearly $L$ is homeomorphic to $\Bbb R$ and is invariant by
$g$.
\end{proof}

\begin{Prop}
Let $G$ be a finitely generated group acting on an $\Bbb R$--tree
$T$ by homeomorphisms. Suppose that the action is non-nesting.
Then if every element of $G$ is elliptic there is an $x\in T$
fixed by $G$.
\end{Prop}
\begin{proof}
We argue by contradiction. Let $G=\langle a_1,a_2,...,a_n \rangle$. If
the intersection $\fix(a_1)\cap \fix(a_2)\cap ... \cap \fix(a_n)= \emptyset $ then
$\fix(a_i)\cap \fix (a_j)=\emptyset $ for some $a_i,a_j$. We claim
that $\smash{a_i^{-1}a_j^{-1}a_ia_j}$ is hyperbolic. Indeed if
$\smash{a_i^{-1}a_j^{-1}a_ia_j(x)=x}$ then $a_ia_j(x)=a_ja_i(x)$. Let
$A=\fix(a_i),B=\fix(a_j)$. We remark that the smallest interval
joining $a_ia_j(x)$ to $A\cup B$ has one endpoint in $A$ while the
smallest interval joining $a_ja_i(x)$ to $A\cup B$ has one
endpoint in $B$ so these two points can not be equal. This is a
contradiction.
\end{proof}
\begin{Prop}
Let $G$ be a  group acting on an $\Bbb R$--tree $T$ by
homeomorphisms. Suppose that the action is non-nesting. Then if
every element of $G$ is elliptic $G$ fixes either an $x\in T$ or
an end of $T$.
\end{Prop}
\begin{proof}
Suppose that $G$ does not fix any $x\in T$. Then there is a
sequence $g_n \in G$ and $x_n\in T$ such that $x_n\in \fix(g_n)$,
$x_n\notin \fix(g_{n-1})$ and $x_n$ goes to infinity. The sequence
$x_n$ defines an end $e$ of $T$. If $r$ is a ray from $x_0\in T$
to $e$ then any $g\in G$ fixes a ray $r_g$ contained in $r$.
Indeed if this is not the case, for some $n$, $\fix(g)$ and
$\fix(g_n)$ are disjoint. It follows as in the previous proposition
that $g^{-1}g_n^{-1}gg_n$ is hyperbolic, which is a contradiction.
\end{proof}

\bibliographystyle{gtart}
\bibliography{link}

\end{document}